\documentclass[a4paper,12pt]{article}

\usepackage[a4paper,hmargin={2.5cm,1.5cm},vmargin={4cm,3cm}]{geometry}

\usepackage[english]{babel}
\usepackage{eurosym}
\usepackage[utf8]{inputenc}
\usepackage{color}
\usepackage{indentfirst}
\usepackage{amsmath}
\usepackage{graphicx}
\usepackage{float}
\usepackage{rotfloat}
\usepackage{rotating}
\usepackage{adjustbox}
\usepackage{multirow}
\usepackage{amssymb}
\usepackage{indentfirst}
 \usepackage{dsfont}
 \usepackage{textcomp}
 \usepackage[T1]{fontenc}
\usepackage{lmodern}
\usepackage{pifont}
\usepackage{graphicx}
\usepackage[pagewise]{lineno}
\newcommand{\R}{\mathbb{R}}
\newtheorem{theorem}{Theorem}[section]

\newtheorem{proposition}[theorem]{Proposition}

\newtheorem{remark}[theorem]{Remark}

 \long\def\symbolfootnote[#1]#2{\begingroup%
\def\thefootnote{\fnsymbol{footnote}}\footnote[#1]{#2}\endgroup}

\begin{document}

\title{%
Improved time-decay for a class of many-magnetic Schr\"odinger flows
}

\author{
Haoran Wang\thanks{%
Department of Mathematics, Beijing Institute of Technology, Beijing 100081, email: wanghaoran@bit.edu.cn}
 } 

\maketitle

\begin{abstract}
\noindent Consider the doubled magnetic Schr\"odinger operator
\begin{equation*}
H_{\alpha,B_0}=\left(i\nabla-\left(\frac{B_0|x|}{2}+\frac{\alpha}{|x|}\right)\left(-\frac{x_2}{|x|},\frac{x_1}{|x|}\right)\right)^2,\quad x=(x_1,x_2)\in\R^2\setminus\{0\},
\end{equation*}
where $\frac{B_0|x|}{2}\left(-\frac{x_2}{|x|},\frac{x_1}{|x|}\right)$ stands for the homogeneous magnetic potential with $B_0>0$ and $\frac{\alpha}{|x|}\left(-\frac{x_2}{|x|},\frac{x_1}{|x|}\right)$ is the well-known Aharonov-Bohm potential with $\alpha\in\R\setminus\mathbb{Z}$.
In this note, we obtain an improved time-decay estimate for the Schr\"odinger flow $e^{-itH_{\alpha,B_0}}$. The key ingredient is the dispersive estimate for $e^{-itH_{\alpha,B_0}}$, which was established in \cite{WZZ23} recently. This work is motivated by L. Fanelli, G. Grillo and H. Kova\v{r}\'{\i}k \cite{FGK15} dealing with the scaling-critical electromagnetic potentials in two and higher dimensions.


\end{abstract}

\baselineskip1.5em

\section{Introduction}

Let us consider the Schr\"odinger operator with electromagnetic potentials
\begin{equation}\label{H-AV}
H_{A, V}=-(\nabla+iA(x))^2+V(x),
\end{equation}
where the electric potential $V(x)$ is a real-valued scalar function on $\R^d$ and the magnetic potential
\begin{equation}
A(x)=(A_1(x),\ldots, A_d(x)): \, \R^d\to \R^d
\end{equation}
is a real-valued vector function on $\R^d$ and fulfils the Coulomb gauge condition
\begin{equation}\label{div0}
\mathrm{div}\, A=0.
\end{equation}
For $d=3$, the vector potential $A$ produces a magnetic field $B$ given by
\begin{equation}\label{B-3}
B=\mathrm{curl} (A)=\nabla\times A.
\end{equation}
For general dimensions $d\geq2$, $B$ should be viewed as a matrix-valued field
$B:\R^d\to \mathcal{M}_{d\times d}(\R)$ given by
\begin{equation}\label{B-n}
B:=DA-DA^t,\quad B_{ij}=\frac{\partial A_i}{\partial x_j}-\frac{\partial A_j}{\partial x_i}.
\end{equation}

The Schr\"odinger operator with electromagnetic potentials \eqref{H-AV} has been extensively studied in spectral and scattering theory; see e.g. the papers of Avron-Herbst-Simon \cite{AHS1,AHS2,AHS3} and the monograph of Reed-Simon \cite{RS}. The study of time-decay estimates for dispersive equations associated to Schr\"odinger operators with electromagnetic potentials has a long history due to their indispensable roles in mathematical physics and partial differential equations(PDEs) ( see e.g. \cite{CS,S}). It turns out that different potentials may result in distinct effects, which means that it is hard to treat all kinds of potentials by a single method. In fact, the picture of thoroughly understanding the electromagnetic Schr\"odinger operators is far from completed, especially those with critical potentials having explicit physical interpretations. For example, the dispersive equations associated to Schr\"odinger operator with Aharonov-Bohm potential (or, the Aharonov-Bohm Hamiltonian) have attracted more and more attentions recently; the Aharonov-Bohm potential is a typical scaling-critical physical model in mathematics in dimension two. In \cite{FFFP13, FFFP15}, Fanelli, Felli, Fontelos and Primo obtained the time-decay estimates for the Schr\"odinger equation associated with the scaling-invariant electromagnetic Schr\"odinger operators including the Aharonov-Bohm Hamiltonian. However, the argument of \cite{FFFP13, FFFP15} fails for the wave equation due to the lack of  pseudoconformal invariance (which was used for Schr\"odinger equation). In \cite{FZZ22}, Strichartz estimate for wave equation was established via the construction of the kernel for the corresponding wave propagator. The authors of \cite{GYZZ22} constructed the spectral measure and further proved the time-decay and Strichartz estimates for Klein-Gordon equation. Nevertheless, the potential models in \cite{FFFP13, FFFP15,FZZ22,GYZZ22} are all scaling-critical ones. Very recently, dispersive estimates for Schr\"odinger equation with two magnetic potentials was established in \cite{WZZ23}; more precisely, the typical 2D magnetic Hamiltonian
\begin{equation}\label{op}
H_{\alpha,B_0}=-(\nabla+i(A_B(x)+A_{\mathrm{hmf}}(x)))^2
\end{equation}
was considered in \cite{WZZ23}, where $A_B(x)$ is the Aharonov-Bohm potential
\begin{equation}\label{AB-potential}
A_B(x)=\alpha\Big(-\frac{x_2}{|x|^2},\frac{x_1}{|x|^2}\Big),\quad x=(x_1,x_2)\in\mathbb{R}^2\setminus\{0\},\alpha\in\mathbb{R}\setminus\mathbb{Z}
\end{equation}
and $A_{\mathrm{hmf}}(x)$ is the homogeneous magnetic potential
\begin{equation}\label{A-hmf}
A_{\mathrm{hmf}}(x)=\frac{B_0}{2}(-x_2,x_1),\quad x=(x_1,x_2)\in\mathbb{R}^2,\quad B_0>0.
\end{equation}
From \eqref{B-n}, we know that the representing function for the magnetic filed in the model \eqref{op} is $B(x)=B_0+\alpha\delta(x)$, where $\delta$ is the usual Dirac delta.
The quadratic form of $H_{\alpha,B_0}$ is positive definite, which means that we are allowed to work with the Friedrichs extension of the operator $H_{\alpha,B_0}$. Associated to this Hamiltonian $H_{\alpha,B_0}$, we can define the Schr\"odinger flow $e^{-itH_{\alpha,B_0}}$ via the Spectral Theorem.
The following dispersive estimate is known to hold for all $t\neq\frac{k\pi}{B_0}$ with $k\in\mathbb{Z}$ (see \cite[Theorem 1.1]{WZZ23})
\begin{equation}\label{dispersive}
\|e^{-itH_{\alpha,B_0}}\|_{L^1(\R^2)\rightarrow L^\infty(\R^2)}\lesssim|\sin(B_0t)|^{-1}.
\end{equation}
In the past few years, many papers were devoted to study the dispersive properties for Schr\"odinger flows associated to Schr\"odinger operators with potentials (mainly electric potentials). For example, dispersive and Strichartz estimates for Schr\"odinger and wave equations associated to Schr\"odinger operator with the inverse-square potential have been obtained in \cite{PST03,BPST03,BPST04} via the techniques such as Morawetz estimates, uniform resolvent estimates, and $TT^\ast$ argument. These arguments fail to work well in the presence of a singular and a unbounded magnetic potentials, even in the pure Aharonov-Bohm case. Nevertheless, several results have been obtained for faster decaying potentials, using the same types of tricks (see e.g. \cite{DF07,EGS08,EGS09,DFVV10,BG12}). In particular, the authors of \cite{GK14} and \cite{FFFP13} independently obtained a representation formula for $e^{-itH_{\alpha,0}}$ in the pure Aharonov-Bohm case and they further proved a polynomial improvement in the decay rate $t^{-1}$, i.e. they obtained the following estimate
\begin{equation}\label{dis:imp}
\||\cdot|^{-\sigma}e^{-itH_{\alpha,0}}|\cdot|^{-\sigma}\|_{L^1(\R^2)\rightarrow L^\infty(\R^2)}\lesssim t^{-1-\sigma},\quad t>0,\quad \sigma\in[0,\mu],
\end{equation}
where $\mu=\mathrm{dist}(\alpha,\mathbb{Z})$. Recently, an explicit representation formula for $e^{-itH_{\alpha,B_0}}$ was obtained in \cite{WZZ23} via two different methods. Motivated by this and the results of \cite{GK14,FGK15} (precisely, \eqref{dis:imp}), we aim to improve the estimate \eqref{dispersive} in this note. More precisely, we will prove the following improved time-decay estimate
\begin{equation}\label{dis:improved}
\||\cdot|^{-\sigma}e^{-itH_{\alpha,B_0}}|\cdot|^{-\sigma}\|_{L^1(\R^2)\rightarrow L^\infty(\R^2)}\lesssim|\sin(B_0t)|^{-1-\sigma},\quad \text{for}\quad t\neq\frac{k\pi}{B_0},
\end{equation}
where $\sigma,\mu$ are the same as \eqref{dis:imp}.
In fact, one can easily verify that $\mu^2$ is the first eigenvalue of the operator $H_{\alpha,0}$ restricted on the unit circle $\mathbb{S}^1$.
We mention that the main ingredient in the proof of \eqref{dis:imp}, apart from the representation formula for the Schr\"odinger kernel, is the dispersive estimate for $e^{-itH_{\alpha,0}}$.

From the mathematical perspective, it is worth emphasizing three key features about the potentials here. First, the Aharonov-Bohm potential is singular at the origin and has the same homogeneity as $\nabla$ so that the perturbation from the Aharonov-Bohm potential \eqref{AB-potential} is non-trivial. Second, the homogeneous magnetic potential \eqref{A-hmf} has degree $1$ so that the Schr\"odinger operator \eqref{op} is no longer scaling invariant. Moreover, the potential $A_{\mathrm{hmf}}(x)$ is unbounded so that the generated magnetic filed $B(x)=B_0$ produces a trapped well. On the one hand, due to the introduction of the magnetic potential \eqref{A-hmf}, the operator $H_{\alpha,B_0}$ has pure point spectrum and the dispersive behavior of the Schr\"odinger equation associated with $H_{\alpha,B_0}$ is quite different from the scaling-invariant models as in \cite{FZZ22, GYZZ22}. On the other, due to the superposition effect from the Aharonov-Bohm potential, a feature of the Mehler kernel, which is related to the Schr\"odinger kernel associated with pure homogeneous magnetic field, breaks down. In fact, when the Aharonov-Bohm effect disappears (i.e. $\alpha=0$), the Schr\"odinger kernel can be written via the classical Mehler formula as
\begin{equation}\label{H0B}
e^{-itH_{0,B_0}}(x,y)=\frac{B_0}{4\pi\sin(B_0t)}\exp\Big\{\frac{B_0}{4i}\big(\cot(B_0t)|x-y|^2-2x\wedge y\big)\Big\},
\end{equation}
where $x\wedge y=x_1y_2-x_2y_1$.
In an equivalent way, \eqref{H0B} can be expressed as
\begin{equation}\label{H0B'}
\begin{split}
e^{itH_{0,B_0}}(x,y)=\frac{B_0}{4\pi\sin(B_0t)}&\exp\Big\{\frac{B_0}{4i}\cot(B_0t)\big(|x|^2+|y|^2\big)\Big\}\\
&\qquad\qquad\times\exp\Big\{i\frac{B_0y\cdot R(B_0t)x}{2\sin(B_0t)}\Big\},
\end{split}
\end{equation}
where $R(\theta)$ is a $2\times 2$ rotation matrix given by
\begin{equation*}
R(\phi)=\left(\begin{array}{cc}\cos\phi &-\sin\phi\\ \sin\phi & \cos\phi \end{array}\right).
\end{equation*}
For our model $H_{\alpha,B_0}$, the kernel of $e^{-itH_{\alpha,B_0}}$ has a similar but different representation formula as \eqref{H0B} (see \eqref{kernel} or \cite[(3.4)]{WZZ23}).


Now we formulate the main result of this note.
\begin{theorem}\label{thm:S}
Let $H_{\alpha,B_0}$ be given by \eqref{op} with magnetic potentials \eqref{AB-potential} and \eqref{A-hmf}.
Let $\mu=\mathrm{dist}(\alpha,\mathbb{Z})$, then, for all $t\neq\frac{k\pi}{B_0},k\in\mathbb{Z}$ and $\sigma\in[0,\mu]$, it holds the following time-decay estimate
\begin{equation}\label{est:dis-improve}
\||\cdot|^{-\sigma}e^{-itH_{\alpha,B_0}}|\cdot|^{-\sigma}\|_{L^1(\R^2)\rightarrow L^\infty(\R^2)}\lesssim|\sin(B_0t)|^{-1-\sigma}.
\end{equation}
In particular, for all $t\in(0,\frac{\pi}{2B_0})$, we have
\begin{equation}\label{dis:smalltime}
\||\cdot|^{-\sigma}e^{-itH_{\alpha,B_0}}|\cdot|^{-\sigma}\|_{L^1(\R^2)\rightarrow L^\infty(\R^2)}\lesssim t^{-1-\sigma}.
\end{equation}
\end{theorem}
\begin{remark}
Theorem \ref{thm:S} implies that a polynomial improvement for the time-decay estimate \eqref{dispersive} is allowed. The optimal decay rate in \eqref{est:dis-improve} depends on the ground level $\mu^2$of the angular part of the Aharonov-Bohm Hamiltonian or the total flux $\alpha$ of the Aharonov-Bohm magnetic field. For this reason, the estimate \eqref{est:dis-improve} is essentially a consequence of the associated dispersive estimate \eqref{dispersive}. On the other hand, since the assumption $\alpha\notin\mathbb{Z}$ implies that $\mu$ is a strictly positive number, \eqref{est:dis-improve} recovers and improves the known time-decay estimate \eqref{dispersive} as $0\leq\sigma\leq\mu$.
In particular, if $\alpha$ is a half-integer (e.g. $\alpha=\frac{1}{2},\frac{3}{2}$, etc.), then we have
\begin{equation*}
\||\cdot|^{-\frac{1}{2}}e^{-itH_{\alpha,B_0}}|\cdot|^{-\frac{1}{2}}\|_{L^1(\R^2)\rightarrow L^\infty(\R^2)}\lesssim|\sin(B_0t)|^{-\frac{3}{2}}.
\end{equation*}
It is interesting to point out that the estimate \eqref{dis:imp} (pure Aharonov-Bohm case) can be viewed as a limit case of \eqref{est:dis-improve} as $B_0\rightarrow0$, in which case the restriction $t\neq\frac{k\pi}{B_0}$ can be removed. The local-in-time estimate \eqref{dis:smalltime} is consistent with \eqref{dis:imp} since it is true that $\frac{2}{\pi}\leq\frac{\sin t}{t}\leq1$ for all $t\in[0,\frac{\pi}{2}]$.
\end{remark}
\begin{remark}
We stress that the restriction $t\neq\frac{k\pi}{B_0}$ cannot be dropped since a trapped well is caused by the unbounded potential \eqref{A-hmf}. Indeed, expanding the square in \eqref{op}, we observe that
\begin{equation*}
\begin{split}
H_{\alpha,B_0}&=-(\nabla+i(A_B(x)+A_{\mathrm{hmf}}(x)))^2
\\&=-\Delta+\frac{B_0^2}{4}|x|^2+\frac{\alpha^2}{|x|^2} +iB_0x^\bot\cdot\nabla+i\frac{2\alpha}{|x|^2}x^\bot\cdot\nabla+\alpha B_0,
\end{split}
\end{equation*}
where $x^\bot=(-x_2,x_1)$.
One will see that the operator is essentially perturbed by the inverse-square potential and the Hermite potential. This phenomenon is natural since the unbounded potential causes a trapped well, the energy cannot be dispersed for long time. This is closely relate to the harmonic oscillator, i.e. $H=-\Delta+|x|^2$, for which Koch and Tataru proved in \cite{KT05}
\begin{equation}
\|e^{-itH}\|_{L^1(\R^d)\to L^\infty(\R^d)}\lesssim|\sin t|^{-\frac d2}.
\end{equation}
For this reason, one observes that the right hand side of \eqref{dis:imp} tends to zero as $t\rightarrow+\infty$, while this is not the case for \eqref{est:dis-improve}. In fact, we take $t=\frac{\pi}{2B_0}+\lambda$ and let $\lambda$ go to infinity, then the right hand side of \eqref{est:dis-improve} is identically equal to $1$ as $t\rightarrow+\infty$.
Finally, we mention that the case $B_0=0$ was considered in \cite[Sect. 6]{Kova12} and a one-dimensional analog of Theorem \ref{thm:S} with inverse square potential was established in \cite{KT14}.
\end{remark}

\section{preliminaries}

In this section, we collect some basic facts about the operator $H_{\alpha,B_0}$ to better understand the magnetic effect.
The space $\mathcal{H}_{\alpha,B_0}^1(\R^2)$ is defined as the completion of $C_c^\infty(\R^2\setminus\{0\};\mathbb{C})$
with respect to the norm
\begin{equation*}
\|\varphi\|_{\mathcal{H}_{\alpha,B_0}^1(\R^2)}=\left(\int_{\R^2}|\nabla_{\alpha,B_0}\varphi(x)|^2\mathrm{d}x\right)^{1/2}
\end{equation*}
where
\begin{equation*}
\nabla_{\alpha,B_0}\varphi=\nabla\varphi+i(A_B+A_{\mathrm{hmf}})\varphi.
\end{equation*}
The quadratic form $Q_{\alpha,B_0}$ of $H_{\alpha,B_0}$ is defined by
\begin{equation*}
\begin{split}
Q_{\alpha,B_0}: & \quad \mathcal{H}_{\alpha,B_0}^1\to \R\\
Q_{\alpha,B_0}(\varphi)&=\int_{\R^2}|\nabla_{\alpha,B_0}\varphi(x)|^2\mathrm{d}x.
\end{split}
\end{equation*}
It is easy to see that $Q_{\alpha,B_0}$ is positive definite, thus the operator $H_{\alpha,B_0}$ is symmetric and semi bounded from below admitting a self-adjoint extension (Friedrichs extension) $H^{F}_{\alpha,B_0}$ with the natural form domain
\begin{equation*}
\mathcal{D}=\Big\{f\in \mathcal{H}_{\alpha,B_0}^1(\R^2):  H^{F}_{\alpha,B_0}f\in L^{2}(\R^2)\Big\}.
\end{equation*}
Although the operator $H_{\alpha,B_0}$ has more than one self-adjoint extension (see \cite{ESV}) by von Neumann's extension theory, we adapt the simplest Friedrichs extension and briefly write the Hamiltonian $H_{\alpha,B_0}$ as its Friedrichs extension $H^{F}_{\alpha,B_0}$.

The spectrum of the Hamiltonian $H_{\alpha,B_0}$ consists of pure discrete eigenvalues and the corresponding ($L^2$-normalized) eigenfunctions form a complete orthonormal basis for $L^2(\R^2)$ (see e.g. \cite{Stov17}).
\begin{proposition}[\cite{Stov17}]
Let $H_{\alpha,B_0}$ be the Hamiltonian given by \eqref{op}.
Then the spectrum of $H_{\alpha,B_0}$ consists of discrete eigenvalues
\begin{equation*}
\lambda_{k,m}=(2m+1+|k+\alpha|+k+\alpha)B_0,\quad k\in\mathbb{Z},\quad m\geq0,
\end{equation*}
and each has a finite multiplicity
\begin{equation*}
\#\Bigg\{j\in\mathbb{Z}:\frac{\lambda_{k,m}-(j+\alpha)B_0}{2B_0}-\frac{|j+\alpha|+1}{2}\in\mathbb{N}\Bigg\}.
\end{equation*}
Furthermore, the corresponding eigenfunctions are given by
\begin{equation*}
V_{k,m}(x)=|x|^{|k+\alpha|}e^{-\frac{B_0 |x|^2}{4}}\, P_{k,m}\Bigg(\frac{B_0|x|^2}{2}\Bigg)e^{ik\theta},\quad \theta=\frac{x}{|x|},\quad k\in\mathbb{Z}, m\geq0,
\end{equation*}
where $P_{k,m}$ is the polynomial of degree $m$ given by
\begin{equation*}
P_{k,m}(r)=\sum_{n=0}^m\frac{(-m)_n}{(1+|k+\alpha|)_n}\frac{r^n}{n!}.
\end{equation*}
with $(a)_n$ ($a\in\R$) denoting the Pochhammer's symbol.
\end{proposition}
\begin{remark}
Recall the generalised Laguerre polynomials $L^\alpha_m(t)$
\begin{equation*}
L^\alpha_m(t)=\sum_{n=0}^m (-1)^n \Bigg(
  \begin{array}{c}
    m+\alpha \\
    m-n \\
  \end{array}
\Bigg)\frac{t^n}{n!}
\end{equation*}
and verify the well-known orthogonality
\begin{equation*}
\int_0^\infty t^{\alpha} e^{-t}L^\alpha_m(t) L^\alpha_n (t)\, \mathrm{d}t=\frac{\Gamma(n+\alpha+1)}{n!} \delta_{n,m},
\end{equation*}
where $\delta_{n,m}$ is the Kronecker delta. One has
\begin{equation}\label{P-L}
P_{k,m}(\rho)=\Bigg(
  \begin{array}{c}
    m+|\alpha+k| \\
    m \\
  \end{array}
\Bigg)^{-1}L^{|\alpha+k|}_m(\rho),
\end{equation}
from which one gets the normalized constant
\begin{equation}\label{V-km-2}
\|V_{k,m}\|^2_{L^2} =\pi \Big(\frac{2}{B_0}\Big)^{1+|\alpha+k|}\Gamma(1+|\alpha+k|) \Bigg(
  \begin{array}{c}
    m+|\alpha+k| \\
    m \\
  \end{array}
\Bigg)^{-1}.
\end{equation}
The Poisson kernel formula for the generalised Laguerre polynomials (see \cite[(6.2.25)]{AAR01})
\begin{equation}\label{La-po}
\begin{split}
&\sum_{m=0}^\infty e^{-cm}\frac{m !}{\Gamma(m+\alpha+1)} L_m^{\alpha} (a) L_m^{\alpha} (b),\qquad a, b, c, \alpha>0\\
&=\frac{e^{\frac{\alpha c}2}}{(ab)^{\frac{\alpha}2}(1-e^{-c})} \exp\left(-\frac{(a+b)e^{-c}}{1-e^{-c}}\right) I_{\alpha}\left(\frac{2\sqrt{ab}e^{-\frac{c}2}}{1-e^{-c}}\right),
\end{split}
\end{equation}
together with \eqref{P-L}, gives
\begin{equation}\label{Po-L}
\begin{split}
&\sum_{m=0}^\infty e^{-cm}\frac{m !}{\Gamma(m+|\alpha+k|+1)}
\Bigg(\begin{array}{c}
    m+|\alpha+k| \\
    m \\
  \end{array}
\Bigg)^2 P_{k,m} (a) P_{k,m} (b)\\
&=\frac{e^{\frac{|\alpha+k|c}2}}{(ab)^{\frac{|\alpha+k|}2}(1-e^{-c})} \exp\left(-\frac{(a+b)e^{-c}}{1-e^{-c}}\right) I_{|\alpha+k|}\left(\frac{2\sqrt{ab}e^{-\frac{c}2}}{1-e^{-c}}\right).
\end{split}
\end{equation}
In view of \eqref{V-km-2}, the formula \eqref{Po-L} actually yields the kernel for the Schr\"odinger flow $e^{-itH_{\alpha,H_0}}$
\begin{equation}\label{kernel}
\begin{split}
e^{-itH_{\alpha,H_0}}(x,y)=&\frac{B_0e^{-itB_0\alpha}}{8\pi^2i\sin(B_0t)}e^{-\frac{B_0(r_1^2+r_2^2)}{4i\tan (tB_0)}}\\
&\quad\times\sum_{k\in\mathbb{Z}}\Bigg( e^{ik(\theta_1-\theta_2-tB_0)}I_{|\alpha+k|}\bigg(\frac{B_0r_1r_2}{2i\sin (tB_0)}\bigg)\Bigg),
\end{split}
\end{equation}
where $x=|x|\frac{x}{|x|}=r_1\theta_1$ and $y=|y|\frac{y}{|y|}=r_2\theta_2$.
\end{remark}

\section{Proof of Theorem \ref{thm:S}}

The main tool for the proof of Theorem \ref{thm:S} is the representation formula
\begin{equation}\label{rep}
\begin{split}
\left(e^{-itH_{\alpha,B_0}}f\right)(x)=&\frac{B_0e^{-itB_0\alpha}}{8\pi^2i\sin(B_0t)}\int_{\R^2}e^{\frac{iB_0(|x|^2+|y|^2)}{4\tan(B_0t)}}\\
&\times\sum_{k\in\mathbb{Z}}\left(e^{ik\left(\frac{x}{|x|}-\frac{y}{|y|}-B_0t\right)}I_{|\alpha+k|}\left(\frac{B_0|xy|}{2i\sin(B_0t)}\right)\right)f(y)\mathrm{d}y
\end{split}
\end{equation}
for any $f\in C_c^\infty(\R^2)$.
In view of the interpolation argument, it is sufficient to verify \eqref{est:dis-improve} for $\sigma=0$ and $\sigma=\mu$. The case $\sigma=0$ is trivially the dispersive estimate \eqref{dispersive}. It remains to verify for the case $\sigma=\mu$.

By the integral representation for the modified Bessel function $I_\nu(z)$ (see \cite[(9.6.18)]{AS65})
\begin{equation}\label{bessel:exp}
I_\nu(z)=\frac{(z/2)^\nu}{\pi\Gamma(\frac{1}{2}+\nu)}\int_{-1}^1(1-s^2)^{\nu-\frac{1}{2}}e^{zs}\mathrm{d}s,\quad z\in\mathbb{C},
\end{equation}
we obtain an upper bound
\begin{equation}\label{bessel:bound}
|I_\nu(iz)|\lesssim\frac{|z|^\nu}{2^\nu\Gamma(\frac{1}{2}+\nu)},\quad \forall z\in\R,\nu\geq0.
\end{equation}
We decompose the whole space $\R^4$ as two parts $\Omega_1,\Omega_2$, where
\begin{equation}\label{decomp-1}
\Omega_1=\left\{(x,y)\in\R^4:|xy|\geq\frac{2\sin(B_0t)}{B_0}\right\}
\end{equation}
and
\begin{equation}\label{decomp-2}
\Omega_2=\left\{(x,y)\in\R^4:0\leq|xy|<\frac{2\sin(B_0t)}{B_0}\right\}.
\end{equation}
By \eqref{dispersive} and \eqref{kernel}, we see that
\begin{equation*}
\sup_{(x,y)\in\R^4}\left|\sum_{k\in\mathbb{Z}}\left(e^{ik\left(\frac{x}{|x|}-\frac{y}{|y|}-tB_0\right)}
I_{|k+\alpha|}\left(\frac{B_0|xy|}{2i\sin(B_0t)}\right)\right)\right|<\infty,\quad \forall t\neq\frac{k\pi}{B_0}.
\end{equation*}
Due to the fact $\mu>0$, we have
\begin{equation*}
K_1:=\sup_{(x,y)\in\Omega_1}\left(\frac{B_0|xy|}{2\sin(B_0t)}\right)^{-\mu}\left
|\sum_{k\in\mathbb{Z}}\left(e^{ik\left(\frac{x}{|x|}-\frac{y}{|y|}-tB_0\right)}
I_{|k+\alpha|}\left(\frac{B_0|xy|}{2i\sin(B_0t)}\right)\right)\right|<\infty
\end{equation*}
On the other hand, we compute
\begin{equation*}
\begin{split}
K_2:&=\sup_{(x,y)\in\Omega_2}\left(\frac{B_0|xy|}{2\sin(B_0t)}\right)^{-\mu}\left
|\sum_{k\in\mathbb{Z}}\left(e^{ik\left(\frac{x}{|x|}-\frac{y}{|y|}-tB_0\right)}
I_{|k+\alpha|}\left(\frac{B_0|xy|}{2i\sin(B_0t)}\right)\right)\right|\\
&\leq\sup_{(x,y)\in\Omega_2}\left(\frac{B_0|xy|}{2\sin(B_0t)}\right)^{-\mu}\sum_{k\in\mathbb{Z}}
\left|I_{|k+\alpha|}\left(\frac{B_0|xy|}{2i\sin(B_0t)}\right)\right|\\
&\lesssim\sup_{(x,y)\in\Omega_2}\sum_{k\in\mathbb{Z}}\frac{1}{2^{|k+\alpha|}\Gamma(\frac{1}{2}+|k+\alpha|)}
\left(\frac{B_0|xy|}{2\sin(B_0t)}\right)^{|k+\alpha|-\mu}.
\end{split}
\end{equation*}
Note that the power of $\frac{B_0|xy|}{2\sin(B_0t)}$ in the series is always positive in view of the definition of $\mu$, we conclude that $K_2<\infty$.

In view of \eqref{decomp-1} and \eqref{decomp-2}, we have obtained
\begin{equation*}
\begin{split}
\sup_{(x,y)\in\R^4}|x|^{-\mu}&\left|\sum_{k\in\mathbb{Z}}\left(e^{ik\left(\frac{x}{|x|}-\frac{y}{|y|}-tB_0\right)}
I_{|k+\alpha|}\left(\frac{B_0|xy|}{2i\sin(B_0t)}\right)\right)\right||y|^{-\mu}\\
&\qquad\leq(B_0/2)^\mu\max\{K_1,K_2\}|\sin(B_0t)|^{-\mu},
\end{split}
\end{equation*}
which yields the estimates \eqref{est:dis-improve} and thus the proof of Theorem \ref{thm:S} is completed.

\vspace{0.2cm}

{\bf Acknowledgement:} The author is currently a graduate student in Beijing Institute of Technology. The author claims that there is no fund supported and no potential conflict of interests.


\vspace{0.5cm}

\begin{thebibliography}{99}

\bibitem{AS65} M. Abramowitz and I. A. Stegun, Handbook of Mathematical Functions with Formulas, Graphs and Mathematical Tables, Dover Publications, New York, 1965.

\bibitem{AHS1} J. Avron, I. Herbst and B. Simon, Schr\"odinger operators with magnetic fields, I. General interactions, Duke Math. J. 45 (1978), 847-883.

\bibitem{AHS2} J. Avron, I. Herbst and B. Simon, Schr\"odinger operators with magnetic fields, II. Separation of center of mass in homogeneous magnetic fields, Ann. Phys. 114 (1978), 431-451.

\bibitem{AHS3} J. Avron, I. Herbst and B. Simon, Schr\"odinger operators in magnetic fields, III. Atoms in homogeneous magnetic field, Commun. Math. Phys. 79 (1981), 529-572.

\bibitem{AAR01} G. E. Andrews, R. Askey and R. Roy, Special Functions (Encyclopedia of Mathematics and its Applications), Cambridge University Press, 2001.

\bibitem{BPST03} N. Burq, F. Planchon, J. Stalker and A. S. Tahvildar-Zadeh, Strichartz estimates for the wave and Schr\"odinger equations with the inverse-square potential, J. Funct. Anal. 203 (2003), no. 2, 519-549.

\bibitem{BPST04} N. Burq, F. Planchon, J. Stalker and A. S. Tahvildar-Zadeh, Strichartz estimates for the wave and Schr\"odinger equations with potentials of critical decay, Indiana Univ. Math. J. 53 (2004), no. 6, 1665-1680.

\bibitem{BG12} M. Beceanu and M. Goldberg, Decay estimates for the Schr\"odinger equation with critical potentials, Comm. Math. Phys. 314 (2012), 471-481.

\bibitem{CS} S. Cuccagna and P. P. Schirmer, On the wave equation with a magnetic potential, Comm. Pure Appl. Math. 54 (2001), 135-152.

\bibitem{DF07} P. D'Ancona and L. Fanelli, Decay estimates for the wave and Dirac equations with a magnetic potential, Comm. Pure Appl. Math. 60 (2007), 357-392.

\bibitem{DFVV10} P. D'Ancona, L. Fanelli, L. Vega, and N. Visciglia, Endpoint Strichartz estimates for the magnetic Schr\"odinger equation, J. Funct. Anal. 258 (2010), 3227-3240.

\bibitem{EGS08} M.B. Erdogan, M. Goldberg and W. Schlag, Strichartz and smoothing estimates for Schr\"odinger operators with large magnetic potentials in $\R^3$, J. Eur. Math. Soc. 10 (2008), 507-531.

\bibitem{EGS09} M.B. Erdogan, M. Goldberg and W. Schlag, Strichartz and smoothing estimates for Schr\"odinger operators with almost critical magnetic potentials in three and higher dimensions, Forum Math. 21 (2009), 687-722.

\bibitem{ESV} P. Exner, P. \v{S}\v{t}ov\'{\i}\v{c}ek and P. Vyt\v{r}as, Generalised boundary conditions for the Aharonov-Bohm effect combined with a homogeneous magnetic field, J. Math. Phys. 43 (2002), 2151.

\bibitem{FFFP13} L. Fanelli, V. Felli, M. A. Fontelos and A. Primo, Time decay of scaling critical electromagnetic Schr\"odinger flows, Comm. Math. Phys. 324 (2013), 1033-1067.

\bibitem{FGK15} L. Fanelli, G. Grillo and H. Kova\v{r}\'{\i}k, Improved time-decay for a class of scaling critical electromagnetic Schr\"odinger flows, J. Funct. Anal. 269 (2015), 3336-3346.

\bibitem{FFFP15} L. Fanelli, V. Felli, M. A. Fontelos and A. Primo, Time decay of scaling invariant electromagnetic Schr\"odinger equations on the plane, Comm. Math. Phys. 337 (2015), 1515-1533.

\bibitem{FZZ22} L. Fanelli , J. Zhang and J. Zheng, Dispersive estimates for 2D-wave equations with critical potentials, Adv. Math. 2022;400.

\bibitem{GK14} G. Grillo and H. Kova\v{r}\'{\i}k, Weighted dispersive estimates for two-dimensional Schr\"odinger operators with Aharonov-Bohm magnetic field, J. Differential Equations 256 (2014), 3889-3911.

\bibitem{GYZZ22} X. Gao , Z. Yin, J. Zhang and J. Zheng, Decay and Strichartz estimates in critical electromagnetic fields, J. Funct. Anal. 2022;282(5).

\bibitem{KT05} H. Koch and D. Tataru, $L^p$ eigenfunction bounds for the Hermite operator, Duke Math. J. 128 (2005), 369-392.

\bibitem{Kova12} H. Kova\v{r}\'{\i}k, Heat kernels of two-dimensional magnetic Schr\"odinger and Pauli operators, Calc. Var. Partial Differential Equations 44 (2012), 351-374.

\bibitem{KT14} H. Kova\v{r}\'{\i}k and F. Truc, Schr\"odinger operators on a half-line with inverse square potentials, Math. Model. Nat. Phenom. 9 (2014), no. 5, 170-176.

\bibitem{PST03} F. Planchon, J. Stalker and A. S. Tahvildar-Zadeh, Dispersive estimates for the wave equation with the inverse-square potential, Discrete Contin. Dyn. Syst. 9 (2003), 1387-1400.

\bibitem{RS} M. Reed and B. Simon, Methods of modern mathematical physics II. Fourier analysis, self-adjointness, Academic Press, New York-London, 1975.

\bibitem{S} W. Schlag, Dispersive estimates for Schr\"odinger operators: a survey. Mathematical aspects of nonlinear dispersive equations, 255-285, Ann. Math. Stud. 163, Princeton Univ. Press, Princeton, NJ, 2007.

\bibitem{Stov17} P. \v{S}\v{t}ov\'{\i}\v{c}ek, The heat kernel for two Aharonov-Bohm solenoids in a uniform magnetic field, Ann. Phys. 376 (2017), 254-282.

\bibitem{WZZ23} H. Wang, F. Zhang and J. Zhang, Decay estimates for one Aharonov-Bohm solenoid in a uniform magnetic field I: Schr\"odinger equation, arXiv:2309.07635v1.


\end{thebibliography}
\end{document}